\RequirePackage[l2tabu, orthodox]{nag}

\documentclass[11pt]{amsart}
\usepackage{fullpage,url,amssymb,enumitem,colonequals,graphicx}
\usepackage[all]{xy} 
\usepackage{mathrsfs} 

\usepackage[dvipsnames,xcdraw,hyperref]{xcolor}

\newcommand{\defi}[1]{\textsf{#1}} 

\newcommand{\C}{\mathbb{C}}
\newcommand{\F}{\mathbb{F}}

\newcommand{\Q}{\mathbb{Q}}

\newcommand{\T}{\mathbb{T}}
\newcommand{\V}{\mathbb{V}}
\newcommand{\Z}{\mathbb{Z}}

\newcommand{\Zhat}{{\hat{\Z}}}

\newcommand{\kbar}{{\overline{k}}}

\newcommand{\Adeles}{\mathbb{A}}

\newcommand{\calO}{\mathcal{O}}

\DeclareMathOperator{\disc}{disc}

\DeclareMathOperator{\End}{End}

\DeclareMathOperator{\Hom}{Hom}

\newcommand{\M}{\operatorname{M}}
\newcommand{\DM}{\operatorname{M}_*}

\newcommand{\Directsum}{\bigoplus} 

\newcommand{\injects}{\hookrightarrow}
\newcommand{\intersect}{\cap} 
\newcommand{\isom}{\simeq}

\newcommand{\tensor}{\otimes} 

\newcommand{\isomto}{\overset{\sim}{\rightarrow}}

\newtheorem{theorem}{Theorem}[section]
\newtheorem{lemma}[theorem]{Lemma}
\newtheorem{corollary}[theorem]{Corollary}

\theoremstyle{definition}
\newtheorem{definition}[theorem]{Definition}

\theoremstyle{remark}
\newtheorem{remark}[theorem]{Remark}

\makeatletter
\g@addto@macro\bfseries{\boldmath} 
\makeatother

\usepackage{microtype}  

\usepackage[
	pagebackref,
	pdfauthor={Bjorn Poonen}, 
]{hyperref}
\usepackage[alphabetic,backrefs,lite,nobysame]{amsrefs} 

\makeatletter
\@namedef{subjclassname@2020}{%
  \textup{2020} Mathematics Subject Classification}
\makeatother

\begin{document}

\title{Lattices in Tate modules}
\author{Bjorn Poonen}
\thanks{This is a corrected version of an article that has published in \emph{Proc.\ Nat.\ Acad.\ Sciences} \textbf{118} (49) e2113201118, \url{https://doi.org/10.1073/pnas.2113201118}\;.
The article arose out of a discussion at the virtual conference ``Arithmetic, Geometry, Cryptography and Coding Theory'' hosted by the Centre International de Rencontres Mathématiques in Luminy in 2021.  B.P.\ was supported in part by National Science Foundation grants DMS-1601946 and DMS-2101040 and Simons Foundation grants \#402472 and \#550033.}
\address{Department of Mathematics, Massachusetts Institute of Technology, Cambridge, MA 02139-4307, USA}
\email{poonen@math.mit.edu}
\urladdr{\url{http://math.mit.edu/~poonen/}}

\author{Sergey Rybakov}
\address{Institute for information transmission problems of the Russian Academy of Sciences}
\address{Interdisciplinary Scientific Center J.-V. Poncelet (ISCP)}
\email{rybakov.sergey@gmail.com}%

\date{April 7, 2025}

\begin{abstract}
Refining a theorem of Zarhin, we prove that given a $g$-dimensional abelian variety $X$ and an endomorphism $u$ of $X$, there exists a matrix $A \in \M_{2g}(\Z)$ such that each Tate module $T_\ell X$ has a $\Z_\ell$-basis on which the action of $u$ is given by $A$, and similarly for the covariant Dieudonn\'e module tensored with $\Q$ if over a perfect field of characteristic $p$.
\end{abstract}

\maketitle

    \section{Introduction}\label{S:introduction}

    Let $X$ be an abelian variety of dimension $g$ over a field $k$ of characteristic $p \ge 0$. 
    Let $\End X$ be its endomorphism ring.
    Let $\End^\circ X \colonequals (\End X) \tensor \Q$.
    Define Tate modules
    \begin{align*}
       T_\ell = T_\ell X &\colonequals \textstyle\varprojlim_n X[\ell^n](\kbar) 
            & V_\ell = V_\ell X &\colonequals T_\ell X \tensor_{\Z_\ell} \Q_\ell \qquad \textup{for each prime $\ell \ne p$} \\
       \T = \T X &\colonequals \prod_{\ell \ne p} T_\ell X 
            & \V = \V X &\colonequals \T X \tensor_\Z \Q \isom {\prod_{\ell \ne p}}' (V_\ell X,T_\ell X);
    \end{align*}
    these are free rank $2g$ modules over $\Z_\ell$, $\Q_\ell$, 
    $\Zhat^{(p)} \colonequals \prod_{\ell \ne p} \Z_\ell$, 
    and $\Adeles^{(p)} \colonequals \Zhat^{(p)} \tensor_\Z \Q \isom \prod^{'}_{\ell \ne p} (\Q_\ell,\Z_\ell)$, 
    respectively.
    If $p>0$ and $k$ is perfect, we have also a $p$-adic analogue of $T_\ell X$, namely the \emph{covariant} Dieudonn\'e module $\DM(X)$, and we define 
    \begin{align*}
       M  &\colonequals \DM(X)
            & M_\Q &\colonequals M \tensor_\Z \Q = M[1/p] \\
       \T_W &\colonequals \T \times M
            & \V_W  &\colonequals \T_W \tensor_\Z \Q = \V \times M_\Q;
    \end{align*}
    these are free rank $2g$ modules over the ring of Witt vectors $W \colonequals W(k)$, its fraction field $K \colonequals W \tensor_\Z \Q =  W[1/p]$, the product $\Zhat^{(p)} \times W$, and $\Adeles_W \colonequals (\Zhat^{(p)}\times W) \tensor_\Z \Q = \Adeles \times \Q$, respectively.
    
    \begin{definition}
    Given rings $R \subseteq R'$ and corresponding modules $L \subseteq L'$, say that $L$ is an \defi{$R$-lattice} in $L'$ if $L$ has an $R$-basis that is an $R'$-basis for $L'$.
    \end{definition}
    
    Zarhin \cite{Zarhin2020}*{Theorem~1.1} proved that given $u \in \End^\circ X$, 
    there exists a matrix $A \in \M_{2g}(\Q)$ such that for every $\ell \ne p$,
    there is a $\Q_\ell$-basis of $V_\ell$ on which the action of $u$ is given by $A$;
    equivalently, there exists a $u$-stable $\Q$-lattice in the $(\prod_{\ell \ne p} \Q_\ell)$-module $\prod_{\ell \ne p} V_\ell$.
    Our main theorem refines this:
    
    \begin{theorem} 
    \label{T:main}
    Let $u \in \End X$.
    \begin{enumerate}[ label=\upshape(\alph*)]
    \item \label{I:V version} There exists a $u$-stable $\Q$-lattice $V \subset \V$.
    \item \label{I:T version} There exists a $u$-stable $\Z$-lattice $T \subset \T$.
    \item \label{I:VW version} If $p>0$ and $k$ is perfect, then there exists a $u$-stable $\Q$-lattice $V \subset \V_W$. 
    \item \label{I:TW version} If $k=\F_p$, then there exists a $u$-stable $\Z$-lattice $T \subset \T_W$.
    \end{enumerate}
    \end{theorem}

    \begin{remark}
    We do not know if \eqref{I:TW version} holds for every perfect field $k$ of characteristic $p>0$.
    \end{remark}
    
    The following answers a question implicit in \cite{Zarhin2020}*{Remark~1.2}:
    
    \begin{corollary}\label{cor1}
    Let $u \in \End X$. Then there exists a matrix $A \in \M_{2g}(\Z)$ such that for every $\ell \ne p$, there is a $\Z_\ell$-basis of $T_\ell X$ on which the action of $u$ is given by $A$,
    and such that if $p>0$ and $k$ is perfect, there is a $K$-basis of $\DM(X) \tensor_{\Z} \Q$ on which the action of $u$ is given by $A$, and such that if $k=\F_p$, there is a $W$-basis of $\DM(X)$ on which the action of $u$ is given by $A$.
    \end{corollary}

    The characteristic~$0$ case of Theorem~\ref{T:main} can be proved by reducing to the case $k=\C$ and taking rational or integral homology \cite{Zarhin2020}*{Remark~1.2}.
    But pairs $(X,u)$ in characteristic $p>0$ cannot always be lifted to characteristic~$0$ \cite{Oort1987}*{Example~14.5}, so the general case does not seem to follow easily from this.
    
    \bigskip
    
    The authors are grateful to a reviewer for suggesting that we prove Theorem~\ref{T:main}(\ref{I:VW version},\ref{I:TW version}), and to Remy van Dobben de Bruyn for pointing out that our argument for Theorem~\ref{T:main}\eqref{I:TW version} works only for $k=\F_p$.
    
    \section{Proof}
    
    \begin{lemma}
    \label{L:semilinear algebra}
    Suppose that $p>0$ and $k$ is perfect.
    Let $L$ be a finite extension of $\Q_p$.
    Let $N$ be an $(L \tensor_{\Q_p} K)$-module with an automorphism $F$ that is $L$-linear and $K$-semilinear with respect to the Frobenius automorphism $\phi$ of $K$.
    Then $N$ is free over $L \tensor_{\Q_p} K$.
    \end{lemma}
    
    \begin{proof}
    The residue field $\ell$ of $L$ is finite, so it has a largest subextension $\ell'$ embeddable in $k$.
    Let $L' \subset L$ be the corresponding unramified extension of $\Q_p$.
    Let $I = \Hom_{\textup{$\Q_p$-algebras}}(L',K)$.
    Then $L' \tensor_{\Q_p} K \isom \prod_{i \in I} K$.
    Applying $L \tensor_{L'}$ yields $L \tensor_{\Q_p} K \isom \prod_{i \in I} L_i$, where each $L_i$ is a \emph{field} since $K$ is absolutely unramified and any tensor product $\ell \tensor_{\ell'} k$ is a field.
    Now $N \isom \Directsum_{i\in I} N_i$, where each $N_i$ is a $L_i$-vector space.
    
    The action of $\phi$ on $K$ induces a permutation $\pi$ of $I$ that is transitive since $L'/\Q_p$ is Galois with group generated by the Frobenius automorphism.
    If $i \in I$ and $j=\pi(i)$, then the compatible actions of $\phi$ of $K$ and $F$ on $N$ induce compatible isomorphisms $L_i \isomto L_j$ and $N_i \isomto N_j$ for each $i$, so $\dim_{L_i} N_i = \dim_{L_j} N_j$.
    It follows by transitivity of $\pi$ that $\dim_{L_i} N_i$ is independent of $i$.
    Thus the module $N = \Directsum_{i\in I} N_i$ is free over $L \tensor_{\Q_p} K \isom \prod_{i \in I} L_i$.
    \end{proof}
    
    \begin{lemma}
    \label{L:E}
    Let $E$ be a number field contained in $\End^\circ X$.
    Let $\calO = E \intersect \End X$.
    Let $h = 2(\dim X)/[E:\Q]$.
    Then
    \begin{enumerate}[label=\upshape(\roman*)]
    \item \label{I:V_ell}
    The $(E \tensor_\Q \Q_\ell)$-module $V_\ell$ is free of rank $h$.
    \item \label{I:M} If $p>0$ and $k$ is perfect, 
    then the $(E \tensor_\Q K)$-module $M_\Q$ is free of rank $h$.
    \item \label{I:T_ell}
    For each $\ell \nmid p \disc \calO$, the $(\calO \tensor_\Z \Z_\ell)$-module $T_\ell$ is free of rank $h$.
    \item \label{I:V}
    The $(E \tensor_\Q \Adeles^{(p)})$-module $\V$ is free of rank $h$.
    \item \label{I:V_W} If $p>0$ and $k$ is perfect, 
    then the $(E \tensor_\Q \Adeles_W)$-module $\V_W$ is free of rank $h$.
    \end{enumerate}
    \end{lemma}
    
    \begin{proof}\hfill
    \begin{enumerate}[label=\upshape(\roman*)]
    \item This is \cite{Ribet1976}*{Theorem~2.1.1}. 
    \item The following proof is essentially a combination of the proofs of \cite{ChCO}*{Proposition~1.4.3.9(1)} and \cite{Ribet1976}*{Theorem~2.1.1}.
    Write $E \tensor_\Q \Q_p \isom \prod_j E_j$ for some finite extensions $E_j$ of $\Q_p$.
    Correspondingly, $M_\Q \isom \Directsum_j M_j$.
    The Frobenius automorphism of the Dieudonn\'e module $M$ induces an $E_j$-linear and $K$-semilinear automorphism of $M_j$.
    By Lemma~\ref{L:semilinear algebra}, the $(E_j \tensor_{\Q_p} K)$-module $M_j$ is free.
    
    On the other hand, by \cite[]{De}*{p.~96, Corollary}, for any $x \in E$, the characteristic polynomials of the actions of $x$ on $M_\Q$ and $V_\ell$ are the same.
    Now repeat the proof of \cite{Ribet1976}*{Theorem~2.1.1}.
    \item 
    Fix $\ell \nmid p \disc \calO$, where $\disc\calO$ is the discriminant of $\calO$.
    For each prime $\lambda$ of $\calO$ dividing $\ell$,
    let $\calO_\lambda \subset E_\lambda$ be the completions of $\calO \subset E$ at $\lambda$. Since $\ell \nmid \disc \calO$, the ring $\calO_\lambda$ is a discrete valuation ring with fraction field $E_\lambda$, and 
    \[
        E \tensor_\Q \Q_\ell \isom \prod_{\lambda | \ell} E_\lambda
        \quad\text{ and }\quad
        \calO \tensor_\Z \Z_\ell \isom \prod_{\lambda|\ell} \calO_\lambda.
    \]
    These induce decompositions 
    \[
        V_\ell = \prod_{\lambda|\ell} V_\lambda
        \quad\text{ and }\quad 
        T_\ell = \prod_{\lambda|\ell} T_\lambda.
    \]
    By \eqref{I:V_ell}, $\dim_{E_\lambda} V_\lambda = h$.
    Since $T_\lambda$ is a torsion-free finitely generated $\calO_\lambda$-module that spans $V_\lambda$,
    it is free of rank $h$ over $\calO_\lambda$.
    Thus $T_\ell = \prod_{\lambda|\ell} T_\lambda$ is free of rank $h$ over $\calO \tensor_\Z \Z_\ell \isom \prod_{\lambda|\ell} \calO_\lambda$.
    \item We have $E \tensor_\Q \Adeles^{(p)} = \prod^{'} (E \tensor_\Q \Q_\ell,\calO \tensor_\Z \Z_\ell)$, so \eqref{I:V} follows from \eqref{I:V_ell} and~\eqref{I:T_ell}.
    \item Similarly, this follows from \eqref{I:V_ell}, \eqref{I:M}, and~\eqref{I:T_ell}.
    \qedhere
    \end{enumerate} 
    \end{proof}
    
    \begin{proof}[Proof of Theorem~\ref{T:main}]\hfill
    \begin{enumerate}[label=\upshape(\alph*)]
    \item We work in the category of abelian varieties over $k$ up to isogeny.
    By \cite{Zarhin2020}*{Theorem~2.4}, $u$ is contained in a subring of $\End^\circ X$ isomorphic to $\prod_i \M_{r_i}(E_i)$ for some number fields $E_i$.
    Then $X$ is isogenous to $\prod Y_i^{r_i}$ for some abelian varieties $Y_i$ with $E_i \subseteq \End^\circ Y_i$.
    If we can find an $E_i$-stable $\Q$-lattice $V_i \subset \V Y_i$ for each $i$, then we may take $V = \prod V_i^{r_i}$.
    In other words, we have reduced to the case that $u \in E \subseteq \End^\circ X$ for some number field $E$.
    By Lemma~\ref{L:E}\eqref{I:V}, 
    \[
        \V = P \tensor_\Q (E \tensor_\Q \Adeles^{(p)})
    \]
    for some $\Q$-vector space $P$.
    Then $V \colonequals P \tensor_\Q E$ is a $u$-stable $\Q$-lattice in $\V$.
    \item
    Given $u \in \End X$, choose $V$ as in~\eqref{I:V version}.
    We have $\Q \intersect \Zhat^{(p)} = \Z[1/p]$, which we interpret as $\Z$ if $p=0$.
    Then $V \intersect \T$ is a $\Z[1/p]$-lattice in $\T$.  Since $\Z[u] \subset \End X$ is a finite $\Z$-module, the $\Z[u]$-submodule generated by any $\Z[1/p]$-basis of $V \intersect \T$ is a $u$-stable $\Z$-lattice.
    \item 
    As in the proof of \eqref{I:V version} we reduce to the case in which
    $u \in E \subseteq \End^\circ X$ for some number field $E$.
    By Lemma~\ref{L:E}\eqref{I:V_W}, 
    \[    \V_W = P \tensor_\Q (E \tensor_\Q \Adeles_W) \] for some $\Q$-vector space $P$. 
    Then $V \colonequals P \tensor_\Q E$ is a $u$-stable $\Q$-lattice in $\V_W$.
    \item Let $V$ be as in \eqref{I:VW version}.  We have $\Q\intersect (\Zhat^{(p)} \times W) =\Z$.
    Then $V \intersect \T_W$ is a $u$-invariant $\Z$-lattice in $\T_W$: to prove that $V \intersect \T_W$ spans $\T_W$ as a $(\Zhat^{(p)} \times W)$-module, we need that $V \intersect M$ spans $M$ as a $W$-module, which holds since $\Q + W = K$ when $k=\F_p$.
    \qedhere
    \end{enumerate}
    \end{proof}

    \section{Generalizations and counterexamples}
    
    In Theorem~\ref{T:main}, suppose that instead of fixing one endomorphism $u$, we consider a \hbox{$\Q$-subalgebra} $R \subset \End^\circ X$
    (or subring $R \subset\End X$) and ask for an $R$-stable $\Q$-lattice (respectively, $\Z$-lattice), i.e., one that is $r$-stable for every $r \in R$.
    \begin{enumerate}[label=\arabic*.]
    \item
    If $R$ is contained in a subring of $\End^\circ X$ isomorphic to $\prod_i \M_{r_i}(E_i)$ for some number fields $E_i$, then the proof of Theorem~\ref{T:main} shows that an $R$-stable lattice exists.
    \item
    Serre observed that if $X$ is an elliptic curve such that $\End^\circ X$ is a quaternion algebra, then for $R = \End^\circ X$, there is no $R$-stable $\Q$-lattice in any $V_\ell$,
    since $R$ cannot act on a $2$-dimensional $\Q$-vector space.
    \item
    If $R$ is assumed to be commutative, then the conclusions of Theorem~\ref{T:main} can still fail.
    For example, suppose that $Y$ is an elliptic curve such that $\End^\circ Y$ is a quaternion algebra $B$, and $X = Y^2$, and 
    \[
        R = \left\{ \begin{pmatrix} a & b \\ 0 & a \end{pmatrix} : a \in \Q \textup{ and } b \in B \right\} \subset \M_2(B) = \End^\circ X.
    \]
    The ideal $\left( \begin{smallmatrix} 0 & B \\ 0 & 0 \end{smallmatrix} \right)$ has square zero,
    so $R$ is commutative.
    For each nonzero $b \in B$, we have 
    \[
        \left( \begin{smallmatrix} 0 & b \\ 0 & 0 \end{smallmatrix} \right) X = 0 \times Y, 
        \quad\text{ so }\quad 
        \left( \begin{smallmatrix} 0 & b \\ 0 & 0 \end{smallmatrix} \right) \V X = 0 \times \V Y,
    \]
    which is of rank~$2$.
    
    Suppose that there is an $R$-stable $\Q$-lattice $V$ in $\V X$.
    Let $U \colonequals V \intersect (0 \times \V Y)$, which is a $\Q$-vector space of dimension at most $2$.
    Then, for every nonzero $b \in B$,
    the image $\left( \begin{smallmatrix} 0 & b \\ 0 & 0 \end{smallmatrix} \right) V$
    is a $2$-dimensional $\Q$-lattice in $0 \times \V Y$, contained in $U$, and hence equal to $U$.
    Thus we obtain a $\Q$-linear injection $$\left( \begin{smallmatrix} 0 & B \\ 0 & 0 \end{smallmatrix} \right) \injects \Hom(V/U,U) \subset \End V.$$
    It is an isomorphism since $$\dim \left( \begin{smallmatrix} 0 & B \\ 0 & 0 \end{smallmatrix} \right) = 4 = \dim \Hom(V/U,U).$$
    Since $\dim_\Q  \left( \begin{smallmatrix} 0 & b \\ 0 & 0 \end{smallmatrix} \right) V = 2$ for each nonzero $b \in B$,
    we have $\dim_\Q f(V) = 2$ for each nonzero $f$ in $\Hom(V/U,U) \subset \End V$, which is absurd.
    Thus there is no $R$-stable $\Q$-lattice in $\V X$.
    \end{enumerate}


\end{document}